\documentclass{amsart}
\usepackage{amscd}

\newtheorem{thm}{Theorem}[section]
\newtheorem{lem}[thm]{Lemma}
\newtheorem{prop}[thm]{Proposition}

\theoremstyle{definition}
\newtheorem{defn}[thm]{Definition}

\theoremstyle{remark}

\numberwithin{equation}{section}

\newcommand{\C}{\mathbb{F}}

\newcommand{\rftz}{\C(z_1,z_2)}
\newcommand{\rfhz}{\C(z_1,z_2,z_3)}

\newcommand{\End}{\operatorname{End}}
\newcommand{\Aut}{\operatorname{Aut}}
\newcommand{\qhat}{\hat{q}}
\newcommand{\g}{\mathfrak{g}}
\newcommand{\fraksl}{\mathfrak{sl}}
\newcommand{\h}{\mathfrak{h}}
\newcommand{\mff}{\mathfrak{f}}
\newcommand{\mfp}{\mathfrak{p}}
\newcommand{\mfb}{\mathfrak{b}}
\newcommand{\mcR}{\mathcal{R}}

\newcommand{\Cc}{\mathcal{C}}
\newcommand{\Cs}{\Cc_\sigma}
\newcommand{\homcs}{\hom_{\Cc_\sigma}}
\newcommand{\Endcs}{\End_{\Cs}}
\newcommand{\tot}{\tilde{\otimes}}

\begin{document}

\title{Generalized Jordanian $R$-matrices of Cremmer-Gervais type}
\author{Robin Endelman}
\address{University of Cincinnati, Cincinnati, OH 45221-0025,
U.S.A.}
\email{endelman@math.uc.edu}
\author{Timothy J. Hodges}
\email{timothy.hodges@uc.edu}
\thanks{The second author was supported in part by NSA grant MDA904-99-1-0026 and by the Charles P. Taft Foundation}

\date{June 9, 2000}

\begin{abstract}
	An explicit quantization is given of certain skew-symmetric solutions of the classical Yang-Baxter equation, yielding a family of $R$-matrices which generalize to higher dimensions the Jordanian $R$-matrices. Three different approaches to their construction are given: as twists of degenerations of the Shibukawa-Ueno Yang-Baxter operators on meromorphic functions; as boundary solutions of the quantum Yang-Baxter equation;  via a vertex-IRF transformation from solutions to the dynamical Yang-Baxter equation.
\end{abstract}
\maketitle

\section*{Introduction}
	Let $\C$ be an algebraically closed field of characteristic zero. The skew-symmetric solutions of the classical Yang-Baxter equation for a simple Lie algebra are classified by the quasi-Frobenius subalgebras; that is, pairs of the form $(\mff, \omega)$ where $\mff$ is a subalgebra and $\omega: \mff\wedge \mff \to \C$ is a nondegenerate 2-cocycle on $\mff$. By a result of Drinfeld \cite{Dr}, the associated Lie bialgebras admit quantizations. This is done by twisting the enveloping algebra $U(\g)[[h]]$ by an appropriate Hopf algebra 2-cocycle. However neither construction lends itself easily to direct calculation and few explicit examples exist to illustrate this theory. The most well-known is the Jordanian quantum group \cite{DM,Ohn} associated to the classical $r$-matrix $E\wedge H$ inside $\mathfrak{sl}(2)\otimes \mathfrak{sl}(2)$. In \cite{GG1}, Gerstenhaber and Giaquinto constructed explicitly the $r$-matrix $r_{\mfp}$ associated to certain maximal parabolic subalgebras $\mfp$ of $\mathfrak{sl}(n)$. In particular for the parabolic subalgebra $\mfp$ generated by $\mfb^+$ and $F_1, \dots F_{n-2}$, their construction yields 
$$
r_\mfp = n \sum_{i<j} \sum_{k=i}^{j-1} E_{k,i}\wedge E_{i+j-k-1,j} + \sum_{i,j}(j-1) E_{j-1,j} \wedge E_{i,i}
$$
In \cite{GG2}, they raise the problem of quantizing this $r$-matrix, in the sense of constructing an invertible $R \in M_n(\C)\otimes M_n(\C)\otimes \C[[h]]$ satisfying the Yang-Baxter equation and of the form $I + h r + O(h^2)$. When $n=2$, the solution is the well-known Jordanian $R$-matrix. Gerstenhaber and Giaquinto construct a quantization of $r_\mfp$ in the $n=3$ case and verify the necessary relations by direct calculation. We give below the quantization of $r_\mfp$ in the general case. Moreover, we are able to give three separate constructions which emphasize the fundamental position occupied by this $R$-matrix. 

	In the first section we construct $R$ (somewhat indirectly) as an extreme degeneration of the Belavin $R$-matrix. We do this by following the construction by Shibukawa and Ueno of solutions of the Yang-Baxter equation for linear operators on meromorphic functions.  In \cite{SU}, they showed that from any solution of Riemann's three-term equation, they could construct such a solution of the Yang-Baxter equation. These solutions occur in three types: elliptic, trigonometric and rational. Felder and Pasquier \cite{FP} showed that in the elliptic case, these operators, after twisting and restricting to suitable finite dimensional subspaces, yield Belavin's $R$-matrices. In the trigonometric case, the same procedure yields the affinization of the Cremmer-Gervais quantum groups; sending the spectral parameter to infinity then yields the Cremmer-Gervais $R$-matrices themselves. Repeating this procedure in the rational case yields the desired quantization of $r_\mfp$, which we shall denote $R_\mfp$.

	In the second section we show that these $R$-matrices occur as boundary solutions of the modified quantum Yang-Baxter equation, in the sense of Gerstenhaber and Giaquinto \cite{GG2}. It was observed in \cite{GG1} that if $\mathfrak{M}$ is the set of solutions  of the modified classical Yang-Baxter equation, then $\mathfrak{M}$ is a locally closed subset of $\mathbb{P}(\g\wedge \g)$ and $\bar{\mathfrak{M}}-\mathfrak{M}$ consists of solutions to the classical Yang-Baxter equation. The element $r_\mfp$ was found to lie on the boundary of the orbit under the adjoint action of $SL(n)$ of the modified Cremmer-Gervais $r$-matrix. In \cite{GG2}, Gerstenhaber and Giaquinto began an investigation into the analogous notion of boundary solutions of the quantum Yang-Baxter equation. They conjectured that the boundary solutions to the classical Yang-Baxter equation described above should admit quantizations which would be on the boundary of the solutions of their modified quantum Yang-Baxter equation. They confirmed this conjecture for the Cremmer-Gervais $r$-matrix in the $\fraksl(3)$ case using some explicit calculations. We prove the conjecture for the general Cremmer-Gervais $r$-matrix by verifying that the matrices $R_\mfp$ do indeed lie on the boundary of the set of solutions to the modified quantum Yang-Baxter equation. 

	In the third section we show that these matrices may also be constructed via a ``Vertex-IRF'' transformation from certain solutions of the dynamical Yang-Baxter equation given in \cite{EV}. This construction is analogous to the original construction of the Cremmer-Gervais $R$-matrices given in \cite{CG}.

	The position of $R_\mfp$ with relation to other fundamental solutions of the YBE and DYBE can be summarized heuristically  by the diagram below.
\begin{equation*}\label{diag}
\begin{CD}
R_B &&&&\qquad\qquad R_F\\
@VVV &&\qquad\qquad@VVV\\
\hat{R}_{CG} @>>> R_{CG}&&\qquad\qquad \hat{R}_{GN} @>>> R_{GN}\\
@VVV	@VVV \qquad\qquad @VVV @VVV\\
R_{B,r} @>>> R_{\mfp} && \qquad\qquad R_{F,r} @>>> R_{GN,r}
\end{CD}
\end{equation*}
$$
\text{YBE} \qquad \qquad \qquad \qquad\qquad \text{DYBE}
$$
On the left hand side, $R_B$ is Belavin's elliptic $R$-matrix; $R_{CG}$ the Cremmer-Gervais $R$-matrix; $\hat{R}_{CG}$ is the affinization of $R_{CG}$ which is also the trigonometric degeneration of the Belavin $R$-matrix; $R_{B,r}$ is a rational degeneration of the Belavin $R$-matrix. The vertical arrows denote degeneration of the coefficient functions (from elliptic to trigonometric and from trigonometric to linear); the horizontal arrows denote the limit as the spectral parameter tends to infinity. On the right hand side, $R_F$ is Felder's elliptic dynamical $R$-matrix; $\hat{R}_{GN}$ and $R_{F,r}$ are trigonometric and rational degenerations; $R_{GN}$ is the Gervais-Neveu dynamical $R$-matrix and $R_{GN,r}$ is a rational degeneration of the Gervais-Neveu matrix given in \cite{EV}. The passage between the two diagrams is performed by Vertex-IRF transformations. The relationships involved in the top two lines of this diagram are well-known \cite{ABB,CG,ES}.  This paper is concerned with elucidating the position of $R_{\mfp}$ in this picture.

	 The authors would like to thank Tony Giaquinto for many helpful conversations concerning  boundary solutions of the Yang-Baxter equation.

\section{Construction of $R_\mfp$}

\subsection{The YBE for operators on function fields}\label{ybeff}

	Recall that if $A$ is an integral domain and $\sigma$ is an automorphism of $A$, then $\sigma$ extends naturally to the field of rational functions $A(x)$ by acting on the coefficients. Denote by $\rftz$ the field of rational functions in the variables $z_1$ and $z_2$. 
Then for any $\sigma \in \Aut \rftz$, and any $i, j \in \{1,2,3\}$, we may define $\sigma_{ij} \in \Aut \rfhz$ by realizing $\rfhz$ as $\C(z_i,z_j)(z_k)$. Set $\Gamma= \Aut \rftz$. Elements  $R = \sum \alpha_i(z_1,z_2)\sigma_i$ of the group algebra $\rftz[\Gamma]$ act as linear operators on $\rftz$ and we may define in this way $R_{ij}$ as linear operators on $\rfhz$.
Thus we may look for solutions of the Yang-Baxter equation $R_{12}R_{13}R_{23} =R_{23}R_{13}R_{12}$ amongst such operators. Denote by $P$ the operator $P\cdot f(z_1,z_2) =f(z_2,z_1)$.

\begin{thm}
	The operator
$$
	R = -\frac{\kappa}{z_1-z_2} P + \left(1 + \frac{\kappa}{z_1-z_2}\right) I
		= I + \frac{\kappa}{z_1-z_2} (I-P)
$$
satisfies the Yang-Baxter equation for any $\kappa \in \C$.
\end{thm}

\begin{proof}Consider an operator of the general form
$$
	R = \alpha(z_1-z_2) P + \beta(z_1-z_2) I
$$
Then it is easily seen that $R$ satisfies the Yang-Baxter equation if and only if
\begin{equation*} 
\alpha(x)\alpha(y) = \alpha(x-y) \alpha(y) + \alpha(x) \alpha(y-x)
\end{equation*}
and
\begin{equation*} 
\alpha(x)\alpha(y)^2 + \beta(y)\beta(-y)\alpha(x+y)
 = \alpha(x)^2 \alpha(y) + \beta(x) \beta(-x) \alpha(x+y)
\end{equation*}
These equations are satisfied when $\alpha(x) = -\kappa/x$ and $\beta(x)=1 - \alpha(x)$. Moreover these are essentially the only such solutions \cite{DH}.
\end{proof}

	In fact, (at least when $\C$ is the field of complex numbers)  this operator is the limit as the spectral parameter tends to infinity of certain solutions of the Yang-Baxter equation with spectral parameter on meromorphic functions constructed by Shibukawa and Ueno. Recall that in \cite{SU}, they showed that operators of the form
$$
	R(\lambda) = G(z_1-z_2, \lambda)P - G(z_1-z_2, \kappa) I
$$
satisfied the Yang-Baxter equation 
$$R_{12}(\lambda_1)R_{13}(\lambda_1+\lambda_2)R_{23}(\lambda_2) =R_{23}(\lambda_2)R_{13}(\lambda_1+\lambda_2)R_{12}(\lambda_1)$$ 
for any $\kappa \in \C$ if $G$ was of the form 
$$
	G(z, \lambda) = \frac{\theta'(0)\theta(\lambda+z)}{\theta(\lambda)\theta(z)}
$$
and $\theta$ satisfied the equation 
\begin{align*}
	\theta(x+y)&\theta(x-y)\theta(z+w)\theta(z-w) + 
	\theta(x+z)\theta(x-z) \theta(w+y)\theta(w-y)\\
	&+ \theta(x+w)\theta(x-w)\theta(y+z)\theta(y-z) = 0
\end{align*}
The principal solution of this equation is $\theta(z) = \theta_1(z)$, the usual theta function (as defined in, say, \cite{WW}), along with the degenerations of the theta functions, $\sin(z)$ and $z$, as one or both of the periods tend to infinity. Felder and Pasquier \cite{FP} showed that in the case where $\theta$ is a true theta function, these operators, when twisted and restricted to suitable subspaces, yield the Belavin $R$-matrices. When $\theta$ is trigonometric, the operator yields in a similar way the affinizations of the Cremmer-Gervais $R$-matrices \cite{DH}. Letting the spectral parameter tend to infinity in a suitable way yields a constant solution of the YBE on the function field which again yields the usual Cremmer-Gervais $R$-matrices on restriction to finite dimensional subspaces. In the rational case, the same twisting and restriction procedure yields the desired quantization of $r_\mfp$.

	When $\theta(z)=z$ we have $G(z,\lambda) = 1/\lambda + 1/z$. Sending $\lambda$ to infinity (and adjusting by a factor of $-\kappa$), we obtain the solution of the Yang-Baxter equation given in the theorem above. Write $R = I+\kappa r$ where $r = (I-P)/(z_1-z_2)$. Then $r$ is a particularly interesting operator. It satisfies the classical Yang-Baxter equation, both forms of the quantum Yang-Baxter equation and has square zero. Its quantization is then just the exponential $\exp \kappa r = I + \kappa r = R$. 

	Let $V_n$ be the space of polynomials in $z_1$ of degree less than $n$. Then we may identify the space $V_n\otimes V_n$ with the subspace of $\rftz$ consisting of polynomials of degree less than $n$ in both  $z_1$ and $z_2$. Since $R\cdot z_1^iz_2^j = z_1^iz_2^j + \kappa(z_1^iz_2^j-z_2^iz_1^j)/(z_1-z_2)$, $R$ restricts to an operator on $V_n\otimes V_n$. With respect to the natural basis,  $R$ has the form
$$
	R(e_i \otimes e_j) = e_i \otimes e_j - \kappa\sum_k \eta(i,j,k) e_k\otimes e_{i+j-k-1}
$$
where$$
\eta(i,j,k) = \begin{cases}
	1 & \text{ if $i\leq k < j$} \\
	-1 & \text{ if $j\leq k < i$} \\
0 & \text{ otherwise }
	\end{cases}
$$

We now apply a simple twist. Define the operator $\tilde{F}_p$ by 
$\tilde{F}_p\cdot f(z_1,z_2)=f(z_1+p, z_2-p)$.

\begin{lem} Let $F=\tilde{F}_p$. Then $F$ and the above $R$ satisfy:
\begin{enumerate}
\item $F_{21}=F_{12}^{-1}$
\item $F_{12}F_{13}F_{23} = F_{23} F_{13}F_{12}$
\item $R_{12}F_{23}F_{13} = F_{13} F_{23}R_{12}$ 
\item $R_{23}F_{12}F_{13} = F_{13} F_{12}R_{23}$
\end{enumerate}
Hence $R_F = F_{21}^{-1}RF_{12}$ also satisfies the Yang-Baxter equation.
\end{lem}

\begin{proof} 
	The four relations are routine verifications. The fact that $R_F$ then satisfies the Yang-Baxter equation is a well-known fact about $R$-matrices extended to this slightly more general situation. 
\end{proof}

Notice that  $F_{21}^{-1}PF_{12}=P$ and $F_{21}^{-1}F_{12}=F^2=\tilde{F}_{2p}$. Taking $p=h/2$ yields 
$$
	R_F = \tilde{F}_{h} +\frac{\kappa}{z_1-z_2+h}(\tilde{F}_{h}-P) 
$$
Notice that
$$
R_F\cdot z_1^iz_2^j = (z_1+h)^i(z_2-h)^j + \kappa \frac{(z_1+h)^i(z_2-h)^j-z_2^iz_1^j}{z_1-z_2+h}
$$
and again $R_F$ restricts to an operator on $V_n\otimes V_n$.

\begin{defn} \label{defrp} Let $n$ be a positive integer. Define
$$
	R_\mfp = \tilde{F}_{h} -\frac{hn}{z_1-z_2+h}(\tilde{F}_{h}-P) 
$$
restricted to $V_n\otimes V_n$.
\end{defn}

Putting all the above together yields the main result.

\begin{thm} For any $h\in \C$ and positive integer $n$, $R_\mfp$ satisfies the Yang-Baxter equation.
\end{thm}

\subsection{Explicit form of $R_\mfp$}

	 We now find an explicit formula for the matrix coefficients of $R_\mfp$ with respect to the natural basis. 

	Define the coefficients of $R_\mfp$ by $R_\mfp \cdot z_1^i z_2^j = \sum_{a,b} R^{ab}_{ij} z_1^a z_2^b$.

\begin{prop}\label{coef}
	The coefficients of $R_\mfp$ are given by
$$
	R^{ab}_{ij} = (-1)^{j-b}\left[\binom{i}{a}\binom{j}{b} 
	+ n\sum_{k}(-1)^{k-a}\binom{i}{k}\binom{j+k-a-1}{b}\eta(j,k,a)\right] h^{i+j-a-b}
$$
\end{prop}

\begin{proof}
Recall that
$$
R_\mfp\cdot z_1^iz_2^j = (z_1+h)^i(z_2-h)^j - hn\frac{(z_1+h)^i(z_2-h)^j-z_2^iz_1^j}{z_1-z_2+h}
$$
For the second term we note that
\begin{multline*}
\frac{z_1^jz_2^i -(z_1+h)^i(z_2-h)^j}{z_1-z_2+h}= \\
\sum_{k,b,a}(-1)^{j+k-a-b}\binom{i}{k}\binom{j+k-a-1}{b}\eta(j,k,a)h^{i+j-a-b-1}z_1^az_2^{b}
\end{multline*}
Combining this with the binomial expansion of the first term yields the assertion.
\end{proof}

The explicit form of this matrix in the case when $n=3$ can be found in \cite[Page 136]{GG2}.

\subsection{The semiclassical limit}

	The operator $R_\mfp$ is a polynomial function of the parameter $h$ of the form $I + r h + O(h^2)$. By working over a suitably extended field, we may assume that $h$ is a formal parameter. Hence $r$ satisfies the classical Yang-Baxter equation. We now verify that $r$ is  the boundary solution $r_\mfp$ associated to the classical Cremmer-Gervais $r$-matrix found by Gerstenhaber and Giaquinto in \cite{GG1}. 

	Recall that their  solution of the CYBE on the boundary of the component containing the modified Cremmer-Gervais $r$-matrix was (up to a scalar)
$$
b_{CG} = n \sum_{i<j} \sum_{k=1}^{j-i} E_{i,j-k+1}\wedge E_{j,i+k} + \sum_{i,j} (n-j)E_{i,i} \wedge E_{j,j+1}.
$$
(Here as usual we are taking the $E_{ij}$ to be the basis of $\End V$ defined by $E_{ij} e_k = \delta_{jk} e_i$ for a fixed basis $\{e_1, \dots, e_n\}$ of $V$; we shall use the convention $x\wedge y = x\otimes y - y \otimes x$). To pass from the $b_{CG}$ to our matrix $r_{\mfp}$, one applies the automorphism $\phi(E_{ij}) = -E_{n+1-j,n+1-i}$. Thus our matrix is again a boundary solution but for a Cremmer-Gervais $r$-matrix associated to a different choice of parabolic subalgebras. 

\begin{thm}
The operator $R_\mfp$ is of the form $I+r_\mfp h +O(h^2)$ where
$$
	r_\mfp \cdot z_1^iz_2^j = n \sum \eta(i,j,k) z_1^k z_2^{i+j-k-1} + i z_1^{i-1} z_2^j - j z_1^i z_2^{j-1}
$$
In particular the matrix representation of $r_\mfp$ with respect to the usual basis is 
$$n \sum_{i<j} \sum_{k=i}^{j-1} E_{k,i}\wedge E_{i+j-k-1,j} + \sum_{i,j}(j-1) E_{j-1,j} \wedge E_{i,i}.
$$
\end{thm}

\begin{proof}
From Proposition \ref{coef}, the coefficients $r_{ij}^{ab}$ are non-zero only when $b=i+j-a-1$ and in this case,
\begin{align*}
	r_{ij}^{a,i+j-a-1} &=\frac{1}{h}R_{ij}^{a,i+j-a-1}=
	(-1)^{a-i+1}\binom{i}{a}\binom{j}{a-i+1} + n\eta(i,j,a)\\
	&= i \delta_{a,i-1} - j \delta _{a, i}+ n\eta(i,j,a).
\end{align*}
Hence
$$
	r_\mfp \cdot z_1^iz_2^j = n \sum \eta(i,j,k) z_1^k z_2^{i+j-k-1} + i z_1^{i-1} z_2^j - j z_1^i z_2^{j-1}.
$$
Thus interpreting $r_\mfp$ as an operator on $V \otimes V$ we get
$$
	r_\mfp \cdot e_i \otimes e_j = n \sum \eta(i,j,k) e_k \otimes e_{i+j-k-1} + (i-1) e_{i-1}\otimes e_j - (j-1) e_i \otimes e_{j-1}.
$$
In matrix form,
$$
	r_\mfp = n \sum_{i<j} \sum_{k=i}^{j-1} E_{k,i}\wedge E_{i+j-k-1,j} + \sum_{i,j}(j-1) E_{j-1,j} \wedge E_{i,i}.
$$
\end{proof}

\section{Boundary solutions of the Yang-Baxter equation}

\subsection{The modified Yang-Baxter equation}
	 In \cite{GG2}, Gerstenhaber and Giaquinto introduced the {\it modified (quantum) Yang-Baxter equation} (MQYBE). An operator $R \in \End V \otimes V$ is said to satisfy the MQYBE if
	$$R_{12}R_{13}R_{23} - R_{23}R_{13}R_{12} = \lambda (P_{123}R_{12} -P_{213}R_{23}) $$
for some nonzero $\lambda$ in  $\C$. Here by $P_{ijk}$ we mean the permutation operator $P_{ijk}(v_1\otimes v_2 \otimes v_3) = v_{\sigma(1)}\otimes v_{\sigma(2)} \otimes v_{\sigma(3)}$ where $\sigma$ is the permutation $(ijk)$. 
	
	Denote by $\mathfrak{R}$ the set of solutions of the YBE in $\End V \otimes V$ and by $\mathfrak{R}'$ the set of solutions of the MQYBE. Then $\mathfrak{R}'$ is a quasi-projective subvariety of $\mathbb{P}(M_{n^2}(\C))$ and $\bar{\mathfrak{R}'} - \mathfrak{R}'$ is contained in $\mathfrak{R}$ \cite{GG2}.  The elements of $\bar{\mathfrak{R}'} - \mathfrak{R}'$ are naturally called {\it boundary} solutions of the YBE. Little is currently known about this set though we conjecture that it contains some interesting $R$-matrices closely related to the quantizations of Belavin-Drinfeld $r$-matrices \cite{ESS}. Let $R$ be a solution of the YBE for which $PR$ satisfies the Hecke equation $(PR-q)(PR+q^{-1}) =0$.   Set $\lambda = (1-q^2)^2/(1+q^2)^2$. Then $Q = (2R + (q^{-1}-q)P)/(q+q^{-1})$ is a unitary solution of the MQYBE. Roughly speaking what we expect to find is the following. If $R$ is a quantization (in the algebraic sense) of a Belavin-Drinfeld $r$-matrix on $\mathfrak{sl}(n)$, then on the boundary of the component of $\mathfrak{R}'$ containing $Q$, we should find the quantization of the skew-symmetric $r$-matrix associated (in the sense of Stolin) with the parabolic subalgebra of $\mathfrak{sl}(n)$ associated to $r$. We prove this conjecture here for the most well-known example, the Cremmer-Gervais $R$-matrices.

	If $R \in \End (V \otimes V)\hat{\otimes}\C[[h]]$ satisfies the QYBE and is of the form $I + h r + O(h^2)$, then $r$ satisfies the classical Yang-Baxter equation and $R$ is said to be a quantization of $r$. The situation for the MQYBE is slightly more complicated and applies only to the $\mathfrak{sl}(n)$ case. Recall that the modified classical Yang-Baxter equation (MCYBE) for an element $r \in \mathfrak{sl}(n)\otimes \mathfrak{sl}(n)$ is the equation
$$
	[r_{12},r_{13}] + [r_{12}, r_{23}] + [r_{13},r_{23}] = \mu \Omega
$$
where $\Omega$ is the unique invariant element of $\wedge^3 \mathfrak{sl}(n)$ (which in the standard representation is the operator $P_{123} - P_{213}$). If $R$ is of the form $I + h r + O(h^2)$ and is a solution of the MQYBE  then $\lambda$ is of the form $\nu h^2 + O(h^3)$ for some scalar $\nu$. If $\nu\neq 0$, then $r$ satisfies the MCYBE. In this case we say that $R$ is a quantization of $r$.

	There is an analogous notion of boundary solution for the classical Yang-Baxter equation. In \cite{GG1}, Gerstenhaber and Giaquinto showed that the matrix $b_{CG}$ lies on the boundary of the component of the set of solutions to the MCYBE containing the modified Cremmer-Gervais classical $r$-matrix. They conjectured that its quantization should lie on the boundary of the component of $\mathfrak{R}'$ containing the modified Cremmer-Gervais $R$-matrix and proved this in the case $n=3$ in \cite{GG2}. We prove now this conjecture in general by showing that $R_\mfp$ lies on the boundary of this component of $\mathfrak{R}'$.

\subsection{The Cremmer-Gervais solution of the MQYBE}

	Consider the linear operator on $\rftz$
$$
	R= \frac{\qhat p z_2}{pz_2-z_1}P +
	 \left( q - \frac{\qhat p z_2}{pz_2-z_1}\right)F_{p}
$$
where $\qhat = q-q^{-1}$ and $F_p\cdot f(z_1,z_2)=f(p^{-1}z_1,p z_2)$.
When restricted to $V_n \otimes V_n$, the above operator becomes the usual 2-parameter Cremmer-Gervais $R$-matrix \cite{DH,H2}. When $p^n=q^2$, this is the original Cremmer-Gervais $R$-matrix which induces a quantization of $SL(n)$ \cite{CG}.

	If $R$ is any solution of the YBE for which $PR$ satisfies the Hecke equation $(PR-q)(PR+q^{-1}) =0$  then $Q = (2R + (q^{-1}-q)P)/(q+q^{-1})$ is a unitary solution of the MQYBE for $\lambda = (1-q^2)^2/(1+q^2)^2$. Hence the operator $Q_{p,q} = (2R -\qhat P)/(q+ q^{-1})$ satisfies the MQYBE. Explicitly,
$$
Q_{p,q} = F_{p} - \frac{\qhat(z_2+ p^{-1}z_1)}{(q+q^{-1})(z_2-p^{-1}z_1)} (F_{p}-P). 
$$
We call the corresponding matrices induced from these operators, the modified Cremmer-Gervais $R$-matrices.

\subsection{Deformation to the boundary}

	Henceforth take $q^2=p^n$. Then the operator $Q_{p,q}$ becomes
$$
Q_p = F_{p} - \frac{(p^n-1)(z_2+ p^{-1}z_1)}{(p^n+1)(z_2-p^{-1}z_1)} (F_{p}-P) 
$$

This is the modified version of the one-parameter Cremmer-Gervais operator described above. Again $Q_p$ may be restricted to the subspace $V_n \otimes V_n$ where its action is given by 
$$
Q_p \cdot z_1^iz_2^j = p^{j-i}z_1^iz_2^j - \frac{(p^n-1)}{(p^n+1)}\sum [\eta(i,j,k) + \eta(i,j,k-1)]p^{j-k}z_1^kz_2^{i+j-k}.
$$
	Fix $h \in \C$ and $p\in \C^*$, define $\tilde{F}_{p,h}$ by $\tilde{F}_{p,h}\cdot f(z_1,z_2) = f(p^{-1}z_1 + p^{-1}h,pz_2-h)$. Define further,

\begin{multline*}
B_{p,h,n} = \tilde{F}_{p,h}- \frac{(p^{n}-1)(pz_2+ z_1)}{(p^{n}+1)(pz_2-z_1-h)}
 (\tilde{F}_{p,h}-P)\\
	+\frac{h(p^{n}-1)(p+1)}{(p^{n}+1)(p-1)(pz_2-z_1-h)} (\tilde{F}_{p,h}-P)\end{multline*}

Note that 
$$B_{1,h,n}= \frac{hn}{(z_2-z_1-h)} (\tilde{F}_{h}-P)+ \tilde{F}_{h}
$$
since $\tilde{F}_h = \tilde{F}_{1,h}$. This is the operator $R_F$ described above (with $\kappa=-hn$) that restricts to $R_\mfp$ on finite dimensional subspaces.

\begin{prop} For all $h$ and $p\not=1$, $B_{p,h,n}$ is a solution of the MQYBE similar to $Q_p$. 
\end{prop}

\begin{proof}
Define a shift operator $\phi_t:\rftz \to \rftz$ by $\phi_t\cdot f(z_1,z_2) = f(z_1-t, z_2-t)$ and let $\phi_t$ act as usual on operators by conjugation. Then, if $F_{p,t} = \phi_t\circ F_p$,
$$
\phi_t \circ Q_p  =  F_{p,t} -
\frac{(p^n-1)(pz_2+ z_1-t(p+1))}{(p^n+1)(pz_2- z_1-t(p-1))} (F_{p,t}-P) $$
Choose $t = h/(p-1)$. Then $\phi_t\circ Q_p =B_{p,h,n}$. This shows that $B_{p,h,n}$ is similar to $Q_p$ and hence satisfies the MQYBE when $p \neq 1$. 
\end{proof}

	Now the restriction of $B_{p,h,n}$ to $V_n \otimes V_n$ is a rational function of $p$ which belongs to $\mathfrak{R}'$ and which for $p=1$ is $R_\mfp$. Thus $R_\mfp$ must be a ``boundary solution'' of the Yang-Baxter equation.

\section{Vertex-IRF transformations and solutions of the dynamical YBE}

	The original construction of the Cremmer-Gervais $R$-matrices was by a generalised kind of change of basis (a ``vertex-IRF transformation") from the Gervais-Neveu solution of the constant dynamical Yang-Baxter equation. Given the above construction of $R_\mfp$ as a rational degeneration of the Cremmer-Gervais matrices, it is natural to expect that $R_\mfp$ should be connected in the same way with some kind of rational degeneration of the Gervais-Neveu matrices. In fact this is precisely what happens. The appropriate solutions to the constant dynamical Yang-Baxter equation (DYBE) were found by Etingof and Varchenko in \cite{EV}. In classifying certain kinds of solutions to the constant DYBE, they found that all such solutions were equivalent to either a generalized form of the Gervais-Neveu matrix or  to a rational version of this matrix. It turns out that $R_\mfp$ is connected via a vertex-IRF transformation with the simplest of this family of rational solutions to the constant DYBE.

	Recall the framework for the dynamical Yang-Baxter equation given in \cite{Hgf}. Let $H$ be a commutative cocommutative Hopf algebra. Let $B$ be an $H$-module algebra with structure map $\sigma: H \otimes B \to B$.
Denote by $\Cc$ the category of right $H$-comodules. Define a new category $\Cs$ whose objects are right $H$-comodules but whose morphisms are
$\homcs (V, W) = \hom_H (V, W \otimes B)$ where $B$ is given a trivial comodule structure. Composition of morphisms is given by the natural embedding of $\hom_H (V, W \otimes B)$ inside $\hom_H (V\otimes B, W \otimes B)$.

	A tensor product $\tot : \Cs \times \Cs \to \Cs$ is defined  on this category in the following way. 
For objects $V$ and $W$, $V\tot W$ is the usual tensor product of $H$ comodules $V \otimes W$. In order to define the tensor product of two morphisms, define first for any $H$-comodule $W$,  a linear twist map $\tau : B \otimes W \to W \otimes B$ by
$$
	\tau(b \otimes w) = w_{(0)} \otimes \sigma(w_{(1)} \otimes b).
$$
where $w \mapsto \sum w_{(0)} \otimes w_{(1)}$ is the structure map of the comodule $W$.
Then for any pair of morphisms $f : V \to V'$ and $g : W \to W'$, define
$$
f \tot g = (1 \otimes m_B)(1 \otimes \tau \otimes 1)(f \otimes g)
$$
Etingof and Varchenko showed in \cite{EV,EV2} that the bifunctor $\tot$  makes $\Cs$ into a tensor category. Let $V \in \Cs$ For any $R \in \Endcs(V \tot V)$ we define elements of $\Endcs (V \tot V \tot V)$, $R_{12} = R \tot 1$ and $R_{23} = 1 \tot R$. Then $R$ is said to satisfy the   $\sigma$-dynamical braid equation ($\sigma$-DBE) if
$R_{12}R_{23}R_{12} = R_{23}R_{12}R_{23}$. If $R$ is a solution of the $\sigma$-DBE then $RP$ satisfies the $\sigma$-dynamical Yang-Baxter equation:
$$
R_{12}R^{12}_{23}R^{123}_{12} = R_{23}R^{23}_{12}R^{132}_{23}
$$
where for instance $R^{132}_{12}=P_{132}R_{12}P_{123}$.

	A vertex-IRF transformation of a solution of the $\sigma$-DBE can then be defined \cite[Section 3.3]{Hgf} as an invertible linear operator $A : V \to V \otimes B$ (that is,  invertible in the sense of the composition of such operators defined above) such that the conjugate operator $R^A = A_2^{-1} A_1^{-1} R A_1 A_2$ is a ``scalar'' operator in the sense that $ R^A(V \otimes V) \subset  V \otimes V\otimes \C$. In this case $R^A$ satisfies the traditional braid equation \cite[Proposition 3.3]{Hgf}. Thus a vertex-IRF transformation transforms a solution of the $\sigma$-DYBE to a solution of the usual YBE.

	Let $T$ be the usual maximal torus of $SL(n)$.   Let $V$ be the standard representation of $SL(n)$ considered as a comodule over $H=\C[T]$ which we may consider as the group algebra of the weight lattice $P$; i.e., $H = \C[K_\lambda \mid \lambda \in P]$. Then $V$ has a basis $\{e_i\}$ of weight vectors with weights $\nu_i$. Denote the structure map by $\rho : V \to V \otimes \C[T]$. Then $\rho(e_i) = e_i \otimes K_{\nu_i}$.
	
	Let $S(\h^*)$ be the symmetric algebra on $\h^*$ and set $B = \text{Frac} (S(\h^*))$. Define an action $\sigma: H \otimes B \to B$ by
$$
	\sigma(K_\lambda \otimes \nu) = \nu -(\lambda, \nu).
$$ 
Denote $\sigma (K_\lambda\otimes b)$ by $b^{\lambda}$. Recall that $(\nu_i, \nu_j) = \delta_{ij}- 1/n$. This fact will be used repeatedly in the calculations below.

	Let $R$ be the matrix $R_\mfp$ defined in Section \ref{ybeff} with $h=1/n$, considered as an operator on the space $V \otimes V$ where $V$ has basis $\{e_1, \dots , e_n\}$. Set $\tilde{R}=RP$ and let $\tilde{R}^{kl}_{ij}$ be the matrix coefficients of $\tilde{R}$ defined by $\tilde{R} \cdot e_i\otimes e_j = \sum_{k,l} \tilde{R}^{kl}_{ij} e_k\otimes e_l$. From Definition \ref{defrp} we have that for any $z_1$ and $z_2$, 
$$
	\sum_{k,l} \tilde{R}^{kl}_{ij} z_1^{k-1} z_2^{l-1} = \alpha(z_1 -z_2)z_1^{i-1} z_2^{j-1} +\beta(z_1-z_2) (z_1+1/n)^{j-1}(z_2-1/n)^{i-1}.
$$
where $\alpha(x)=1/(x+1/n)$ and $\beta(x) = 1 -\alpha(x)$.
Define the operator $\mcR \in \Endcs V \tot V$ by
\begin{align*}\mcR(e_i \otimes e_j) &=  e_i \otimes e_j \otimes  \alpha(\nu_i^{\nu_j} -\nu_j)
	+e_j \otimes e_i \otimes \beta(\nu_i^{\nu_j} -\nu_j) \\
	&=e_i \otimes e_j \otimes  \frac{1}{\nu_i -\nu_j+\delta_{ij}}
	+e_j \otimes e_i \otimes \left(1 - \frac{1}{\nu_i -\nu_j+\delta_{ij}}\right).
\end{align*}
This is the solution of the DBE corresponding to the standard example of solution of the DYBE of the type given in \cite[Theorem 1.2]{EV}. Finally define an operator $A \in \Endcs(V)$ by
$A(e_i) = \sum e_k \otimes \nu_k^{i-1}$.

\begin{thm}
	$\mcR^A = \tilde{R}$
\end{thm}

\begin{proof} We prove that $\mcR A_1 A_2 = A_1 A_2 \tilde{R}$. In matrix form this is equivalent to
$$
\sum_{c,d} \mcR^{ms}_{cd} (A^c_i)^{\nu_d} A_j^d = \sum_{k,l} \tilde{R}^{kl}_{ij} (A^m_k)^{\nu_s} A_l^s.
$$
Using the fact that $\beta(\nu_m^{\nu_s} - \nu_s) = 0$ when $m=s$ 
\begin{align*}
\sum_{k,l} &\tilde{R}^{kl}_{ij} (A^m_k)^{\nu_s} A_l^s  = \sum_{k,l}\tilde{R}^{kl}_{ij}(\nu_m^{\nu_s})^{k-1}\nu_s^{l-1}\\
	&= \alpha(\nu_m^{\nu_s} - \nu_s)(\nu_m^{\nu_s})^{i-1}\nu_s^{j-1} +\beta(\nu_m^{\nu_s} - \nu_s)(\nu_s-\frac{1}{n})^{i-1}(\nu_m^{\nu_s}+\frac{1}{n})^{j-1}\\
	&= \alpha(\nu_m^{\nu_s} - \nu_s)(\nu_m^{\nu_s})^{i-1}\nu_s^{j-1} +\beta(\nu_m^{\nu_s} - \nu_s)(\nu_s^{\nu_m})^{i-1}(\nu_m)^{j-1}\\
	&= \sum_{c,d} \mcR^{ms}_{cd} (A^c_i)^{\nu_d} A_j^d
\end{align*}
as required.
\end{proof}

\end{document}